# Average Error of the Prime Number Theorem

Li An-Ping

Beijing 100085, P.R. China
apli0001@sina.com

Abstract

In this paper, we will give some estimation about the average error of the prime number theorem.

## 1. Introduction

The functions $\pi(x)$, $\vartheta(x)$, $\psi(x)$, $\Lambda(n)$ are the familiar ones in the number theory

$$\pi(x) = \sum_{p \le x} 1,$$

$$\Lambda(n) = \begin{cases} \log p, & \text{if } n \text{ is a power of prime number } p. \\ 0, & \text{otherwise.} \end{cases}$$

$$\psi(x) = \sum_{n \le x} \Lambda(n),$$

$$\vartheta(x) = \sum_{p \le x} \log p.$$

Where variable $p$ in the sums are of prime numbers. The prime number theorems are the ones to describe the distribution rule of prime numbers, there are mainly the following three equivalent versions

$$\begin{aligned} &\text{i)} \quad \pi(x) \sim x/\log x, \quad (\text{or, Li } x) \\ &\text{ii)} \quad \vartheta(x) \sim x, \\ &\text{iii)} \quad \psi(x) \sim x. \end{aligned} \tag{1.1}$$

That is, $x/\log x$, $x$, $x$ are the dominant parts of functions $\pi(x)$, $\vartheta(x)$ and $\psi(x)$ respectively, when $x$ is sufficient great. So, the rest work is to determine the error parts. In the paper, we are just concerning to the version iii).

Let $r(x) = \psi(x) - x$, de la Vallée Poussin proved that

$$r(x) \ll x \exp\left(-c\sqrt{\log x}\right). \tag{1.2}$$

A. Walfisz improved that

$$r(x) \ll x(\log x)^{6/5} \exp\left(-c_1(\log x)^{3/5}(\log\log x)^{-1/5}\right). \tag{1.3}$$

Where $c$, $c_1$ are two constants.

With $\Omega$-theorems, it is known that the signs of $r(x)$ are changed frequently, so it will be meaningful to investigate the average error $\overline{r}(x) = \frac{1}{x}\sum_{n \le x} r(n)$. In this paper, we will give some estimation for $\overline{r}(x)$. The argument is essentially an extension of Perron's formula.

**Theorem 1.** Suppose that $\{\rho\}$ is the set of non-trivial zeros of Riemann's zeta-function $\zeta(s)$, $\eta = \max\left\{\operatorname{Re}\rho \,|\, |\operatorname{Im}\rho| \le x^{1/4}\log x\right\}$, then

$$\bar{r}(x) \ll x^{\eta}. \tag{1.4}$$

In general,

$$\bar{r}(x) \ll x\exp\left(-2^{2/5}\cdot c_1(\log x)^{3/5}(\log\log x)^{-1/5}\right). \tag{1.5}$$

It should be mentioned that though the factor $(\log x)^{6/5}$ in (1.3) actually may be assimilated into $\exp()$ of its right-hand, in order to coincide with the constant $c_1$ in (1.3) and (1.5), we have preserved it.

**2. The Proof of Theorems 1.**

**Lemma 1.** For $b, T > 0$, it has

$$\frac{1}{2\pi i}\int_{b-iT}^{b+iT}\frac{a^s}{s(s+1)}ds = \begin{cases} 1-1/a+R_1, & \text{if } a>1 \\ R_2, & \text{if } 0<a<1, \\ 1/\pi T + R_3, & \text{if } a=1. \end{cases} \tag{2.1}$$

Where

$$R_1 = O(a^b\min(1/T, 1/(T^2\log a)),\ R_2 = O(a^b\min(1/T, 1/(T^2\log a)),\ R_3 = O(1/T^3).$$

Proof. For $a>1$, let $\mathscr{L}_1$ be the contour of rectangle with vertices $b\pm iT, -U\pm iT$, then by residue theorem, it has

$$\frac{1}{2\pi i}\int_{\mathscr{L}_1}\frac{a^s}{s(s+1)}ds = 1-1/a.$$

Besides,

$$\left|\int_{-U+it}^{b+it}\frac{a^s}{s(s+1)}ds\right| \le \int_{-U}^{b}\frac{a^{\sigma}}{\sqrt{\sigma^2+T^2}\sqrt{(1+\sigma)^2+T^2}}d\sigma \le \frac{1}{T^2}\int_{-U}^{b}a^{\sigma}d\sigma \le \frac{a^b}{T^2\log a}$$

$$\left|\int_{-U-iT}^{-U+iT}\frac{a^s}{s(s+1)}ds\right| \le a^{-U}\int_{-T}^{T}\frac{1}{(U-1)^2+t^2}dt = \frac{a^{-U}}{U-1}\arctan\frac{T}{U-1} = O\left(\frac{Ta^{-U}}{U^2}\right)$$

Moreover, take the contour $\mathscr{L}_2$ with the line segment from $b-iT$ to $b+iT$ and the arc $C$ of its left-hand with the radius $\sqrt{b^2+T^2}$ and the center at the origin. Then by residue theorem,

$$\frac{1}{2\pi i}\int_{\mathscr{L}_2}\frac{a^s}{s(s+1)}ds = 1-1/a,$$

$$\left|\frac{1}{2\pi i}\int_C \frac{a^s}{s(s+1)}ds\right| \le \frac{a^b}{T}.$$

The treatment for the case $0<a<1$ is similar，let $\mathscr{L}_3$ be the contour of rectangle with vertices $b \pm iT, U \pm iT$, by residue theorem it has

$$\frac{1}{2\pi i}\int_{\mathscr{L}_3} \frac{a^s}{s(s+1)}ds = 0,$$

and

$$\left|\int_{U+it}^{b+it} \frac{a^s}{s(s+1)}ds\right| \le \int_U^b \frac{a^\sigma}{\sqrt{\sigma^2+T^2}\sqrt{(1+\sigma)^2+T^2}}d\sigma \le \frac{1}{T^2}\int_U^b a^\sigma d\sigma \le \frac{a^b}{T^2\log a}$$

$$\left|\int_{U-iT}^{U+iT} \frac{a^s}{s(s+1)}ds\right| \le a^U \int_{-T}^{T} \frac{1}{(U+1)^2+t^2}dt = \frac{a^U}{U+1}\arctan\frac{T}{U+1} = O\left(\frac{a^U T}{U^2}\right)$$

Again, take the contour $\mathscr{L}_4$ with line segment from $b-iT$ to $b+iT$ and the arc $C$ of the right-hand with the radius $\sqrt{b^2+T^2}$ the center at the origin. Then,

$$\frac{1}{2\pi i}\int_{\mathscr{L}_4} \frac{a^s}{s(s+1)}ds = 0, \qquad \left|\frac{1}{2\pi i}\int_C \frac{a^s}{s(s+1)}ds\right| \le \frac{a^b}{T}.$$

For $a=1$,

$$\frac{1}{2\pi i}\int_{b-iT}^{b+iT} \frac{1}{s(s+1)}ds = \frac{1}{2\pi}\int_{-T}^{T}\left(\frac{b-it}{b^2+t^2} - \frac{b+1-it}{(b+1)^2+t^2}\right)dt$$

$$= \frac{1}{\pi}\int_0^T\left(\frac{b}{b^2+t^2} - \frac{b+1}{(b+1)^2+t^2}\right)dt = \frac{1}{\pi}\left(\arctan\frac{T}{b} - \arctan\frac{T}{b+1}\right)$$

$$= \frac{1}{\pi T} + O\left(1/T^3\right)$$

The proof is completed. □

Suppose that $A(s)=\sum_{n=1}^{\infty} a(n)n^{-s},\ |a(n)| \le \tau(n),\ \sum_{n=1}^{\infty}|a(n)|n^{-\sigma} \le \Theta(\sigma),\ \sigma \ge \sigma_a$, where $\tau(n)$, $\Theta(\sigma)$ are two real functions, $\tau(n)$ is non-decreasing, $\sigma_a$ is the bound of absolute convergence of $A(s)$. Denote by $\mathcal{A}(x,s)=\sum_{n\le x} a(n)n^{-s},\ F(x,s)=\sum_{n\le x}\mathcal{A}(n,s)$, write $\overline{x}=[x+1]$.

**Lemma 2.** For arbitrary $s_0$, $s_0=\sigma_0+it_0$, and $b$, $b+\sigma_0>\sigma_a$, there is

$$F(x,s_0)=\frac{\overline{x}}{2\pi i}\int_{b-iT}^{b+iT}A(s+s_0)\frac{\overline{x}^s}{s(s+1)}ds+Q. \qquad (2.2)$$

$$Q=O\left(\frac{x^{1+b}\Theta(b+\sigma_0)}{T^2}\right)+O\left(\frac{x^{2-\sigma_0}}{T}\tau(2\overline{x})\min\left(1,\frac{\log x}{T}\right)\right)+O\left(x^{1-\sigma_0}\tau(\overline{x})\frac{1}{T}\right).$$

Proof. By Lemma 1,

$$\frac{1}{2\pi i}\int_{b-iT}^{b+iT}A(s+s_0)\frac{\overline{x}^s}{s(s+1)}ds=\sum_{n=1}^{\infty}a(n)n^{s_0}\frac{1}{2\pi i}\int_{b-iT}^{b+iT}\frac{1}{s(s+1)}\left(\frac{\overline{x}}{n}\right)^s ds$$

$$=\sum_{n\le x}a(n)n^{s_0}(1-n/\overline{x})+R.$$

$$R=\sum_{n=1}^{\infty}|a(n)|n^{\sigma_0}\left(\frac{\overline{x}}{n}\right)^b\frac{1}{T}\min\left(1,\frac{1}{T\log(\overline{x}/n)}\right)$$

$$=\sum_{n\le\overline{x}/2}+\sum_{\overline{x}/2<n\le 2\overline{x}}+\sum_{n>2x}$$

$$\sum_{n\le x/2}=\sum_{n\le x/2}|a(n)|n^{-\sigma_0}\left(\frac{x}{n}\right)^b\frac{1}{T}\min\left(1,\frac{1}{T\log(x/n)}\right)\ll\frac{x^b}{T^2}\sum_{n\le x/2}|a(n)|n^{-\sigma_0-b}$$

$$\sum_{n>x/2}=\sum_{n>x/2}|a(n)|n^{-\sigma_0}\left(\frac{x}{n}\right)^b\frac{1}{T}\min\left(1,\frac{1}{T\log(x/n)}\right)\ll\frac{x^b}{T^2}\sum_{n\le x/2}|a(n)|n^{-\sigma_0-b}$$

So,

$$\sum_{n\le x/2}+\sum_{n>2x}\ll\frac{x^b}{T^2}\Theta(b+\sigma_0).$$

Moreover,

$$\sum_{x/2<n\le 3x/2}\ll x^{-\sigma_0}\tau(\overline{x})\frac{1}{T}+x^{-\sigma_0}\tau(3\overline{x}/2)\min\left(x,\sum_{\substack{x/2<n\le 3x/2\\ n\ne\overline{x}}}\frac{1}{T\log(\overline{x}/n)}\right)$$

As $\left|\log(1+x)\right|=\left|\int_0^x\frac{dz}{1+z}\right|\ge\frac{x}{1+x}$, it follows

$$\sum_{x/2\le n<\overline{x}}\frac{1}{\log(\overline{x}/n)}\le\sum_{x/2\le n\le\overline{x}}\frac{1}{\log(1+(\overline{x}-n)/n)}\le\sum_{1\le r\le\overline{x}/2}\frac{1}{\log(1+r/(\overline{x}-r))}$$

$$\le\overline{x}\sum_{1\le r\le\overline{x}/2}\frac{1}{r}\ll x\log x$$

Similarly,

$$\sum_{\overline{x}<n<3\overline{x}/2}\frac{1}{\log(x/n)}\ll x\log x.$$

Hence,

$$\frac{\bar{x}}{2\pi i}\int_{b-iT}^{b+iT} A(s+s_0)\frac{\bar{x}^s}{s(s+1)}ds=\sum_{n\le x}a(n)n^{-s_0}(\bar{x}-n)+\bar{x}R$$

$$=\sum_{n\le x}\sum_{m\le n}a(m)m^{-s_0}+\bar{x}R$$

$$=\sum_{n\le x}\mathcal{A}(n,s_0)+O\left(\frac{x^{1+b}}{T^2}\Theta(b+\sigma_0)\right)+O\left(\frac{x^{2-\sigma_0}}{T}\tau(2\bar{x})\min\left(1,\frac{\log x}{T}\right)\right)+O\left(\frac{x^{1-\sigma_0}\tau(\bar{x})}{T}\right)$$

The proof is finished. □

The Proof of Theorem 1：

By Euler product formula of Riemann's zeta-function $\zeta(s)$, it has

$$\frac{\zeta'(s)}{\zeta(s)}=\frac{d}{ds}\log\zeta(s)=-\sum_p\frac{p^{-s}\log p}{1-p^{-s}}=-\sum_p\sum_{m\ge1}\log p/p^{ms}=-\sum_n\Lambda(n)/n^s.$$

Denote by $\Psi(x)=\sum_{n\le x}\psi(n)$, and let $b=1+1/\log x$, by Lemma 2, there is

$$\Psi(x)=\frac{\bar{x}}{2\pi i}\int_{b-iT}^{b+iT}\frac{-\zeta'(s)}{\zeta(s)}\frac{\bar{x}^s}{s(s+1)}ds+Q. \qquad (2.3)$$

Where $Q=O\left(\frac{x^2\log T}{T^2}\right)+O\left(\frac{x^2\log^2 x}{T^2}\right)+O\left(\frac{x\log x}{T}\right)$.

Let $\mathscr{L}$ be the contour of rectangle with vertices $b\pm iT,-U\pm iT$. By residue theorem, it has

$$\frac{1}{2\pi i}\int_{\mathscr{L}}\frac{-\zeta'(s)}{\zeta(s)}\frac{\bar{x}^s}{s(s+1)}ds=\frac{\bar{x}}{2}-\frac{\zeta'(0)}{\zeta(0)}+\frac{\zeta'(-1)}{\zeta(-1)}\bar{x}^{-1}-\sum_{n\le U/2}\frac{\bar{x}^{-2n}}{2n(2n-1)}-\sum_{|\gamma|\le T}\frac{\bar{x}^\rho}{\rho(\rho+1)}$$

$$=\frac{\bar{x}}{2}-\sum_{|\gamma|\le T}\frac{\bar{x}^\rho}{\rho(\rho+1)}+O(1).$$

i.e.

$$\frac{1}{2\pi i}\int_{b-iT}^{b+iT}\frac{-\zeta'(s)}{\zeta(s)}\frac{\bar{x}^s}{s(s+1)}ds=\frac{1}{2\pi i}\left(\int_{-U-iT}^{-U+iT}+\int_{b-iT}^{-U-iT}+\int_{-U+iT}^{b+iT}\right)\frac{-\zeta'(s)}{\zeta(s)}\frac{\bar{x}^s}{s(s+1)}ds$$

$$+\frac{\bar{x}}{2}-\sum_{|\gamma|\le T}\frac{\bar{x}^\rho}{\rho(\rho+1)}+O(1)$$

On the other hand,

$$\int_{-U\pm iT}^{b\pm iT}\frac{-\zeta'(s)}{\zeta(s)}\frac{\bar{x}^s}{s(s+1)}ds\ll\frac{(\log T)^2}{T^2}\int_{-U}^{b}\bar{x}^\sigma d\sigma\ll\frac{x^b(\log T)^2}{T^2\log x}$$

$$\int_{-U-iT}^{-U+iT} \frac{-\zeta'(s)}{\zeta(s)} \frac{\overline{x}^s}{s(s+1)} ds \ll (\log T)^2 \overline{x}^{-U} \int_{-T}^{T} \frac{dt}{U^2+t^2} \ll (\log T)^2 / x^{U+1}$$

With (2.3), there is

$$\Psi(x) = \frac{\overline{x}^2}{2} + Q + O\left(\frac{x^2(\log T)^2}{T^2 \log x}\right) + O\left(x^{-U}(\log T)^2\right) - \sum_{|\mathrm{Im}\,\rho| \le T} \frac{\overline{x}^{1+\rho}}{\rho(\rho+1)}$$

i.e.

$$\overline{r}(x) = Qx^{-1} + O\left(\frac{x(\log T)^2}{T^2 \log x}\right) - \sum_{|\mathrm{Im}\,\rho| \le T} \frac{\overline{x}^{\rho}}{\rho(\rho+1)}$$

Suppose that $\eta = \max\left\{\mathrm{Re}\,\rho \,\middle|\, |\mathrm{Im}\,\rho| \le T\right\}$, then

$$\overline{r}(x) = Qx^{-1} + O\left(\frac{x(\log T)^2}{T^2 \log x}\right) + x^{\eta} \sum_{|\mathrm{Im}\,\rho| \le T} \frac{1}{|\rho(\rho+1)|}$$

Take $T = x^{1/4} \log x$, and notice that

$$\sum_{|\mathrm{Im}\,\rho| \le T} \frac{1}{|\rho(\rho+1)|} \le \sum_{|\gamma| \le T} \frac{1}{\gamma^2} \ll \int_1^{\infty} \frac{\log(t)}{t^2} dt = 1.$$

So, (1.4) is followed.

Moreover, A. Walfisz [5] presented that there is a constant $c_0$ such that $\zeta(s)$ is free of zeros in the region

$$\sigma \ge 1 - c_0(\log^2 t \log\log t)^{-1/3}.$$

This means

$$\eta < 1 - c_0(\log^2 T \log\log T)^{-1/3}.$$

Take $\log T = \frac{c_0^{3/5}}{2^{3/5}} \left(\frac{3}{5}\right)^{-1/5} \left(\log x\right)^{3/5} (\log\log x)^{-1/5}$, then

$$\overline{r}(x) \ll x \exp\left(-2^{2/5} c_0^{3/5} \left(3/5\right)^{-1/5} \left(\log x\right)^{3/5} (\log\log x)^{-1/5}\right)$$
$$\ll x \exp\left(-2^{2/5} c_1 \left(\log x\right)^{3/5} (\log\log x)^{-1/5}\right), \qquad c_1 = c_0^{3/5} \left(3/5\right)^{-1/5}.$$

And (1.5) is proved. □

Finally, we wish to mention that the concept of average error may be developed to multilayer

$$\overline{r}^{(k)}(x) = \frac{1}{C_{[x]+k-1}^k} \sum_{n_1 \le x} \sum_{n_2 \le n_1} \cdots \sum_{n_k \le n_{k-1}} r(n_k). \tag{2.6}$$

And like Lemmas 1, 2, consider complex integral

$$\frac{1}{2\pi i}\int_{b-iT}^{b+iT}\frac{a^s}{s\cdot C_{s+k}^{k}}ds.$$

For example, $k=2,3,$

$$\overline{r}^{(2)}(x)=\frac{1}{C_{[x]+1}^{2}}\sum_{n_1\le x}\sum_{n_2\le n_1} r(n_2), \qquad \overline{r}^{(3)}(x)=\frac{1}{C_{[x]+2}^{3}}\sum_{n_1\le x}\sum_{n_2\le n_1}\sum_{n_3\le n_2} r(n_3).$$

Instead of Lemma1 to consider complex integrals respectively

$$\frac{2!}{2\pi i}\int_{b-iT}^{b+iT}\frac{a^s}{s(s+1)(s+2)}ds \quad \text{and} \quad \frac{3!}{2\pi i}\int_{b-iT}^{b+iT}\frac{a^s}{s(s+1)(s+2)(s+3)}ds.$$

With a similar argument, it can be followed that

$$\begin{aligned}&\overline{r}^{(2)}(x)\ll x\exp\left(-3^{2/5}c_1(\log x)^{3/5}(\log\log x)^{-1/5}\right),\\ &\overline{r}^{(3)}(x)\ll x\exp\left(-4^{2/5}c_1(\log x)^{3/5}(\log\log x)^{-1/5}\right).\end{aligned} \tag{2.7}$$

Besides, if Riemann's hypothesis is true, then it has

$$\overline{r}^{(i)}(x)\ll \lambda_i x^{1/2}, \qquad i=1,2,3, \tag{2.8}$$

where

$$\lambda_1=\sum_{|\gamma|\le T}\frac{x^{\gamma i}}{\rho(\rho+1)}, \quad \lambda_2=\sum_{|\gamma|\le T}\frac{x^{\gamma i}}{\rho(\rho+1)(\rho+2)}, \quad \lambda_3=\sum_{|\gamma|\le T}\frac{x^{\gamma i}}{\rho(\rho+1)(\rho+2)(\rho+3)}.$$

In this case, the decisive problem is to find the values $\lambda_i$, $i=1,2,3.$

On the other hand, it is easy to follow that

$$\begin{aligned}\overline{r}^{(i)}(n)&=\frac{1}{C_{n+i-1}^{i}}\sum_{m\le n}C_{n+i-m-1}^{i-1}r(m)=\frac{1}{C_{n+i-1}^{i}}\sum_{m\le n}C_{n+i-m-1}^{i-1}\sum_{j\le m}(\Lambda(j)-1)\\ &=\frac{1}{C_{n+i-1}^{i}}\left(\sum_{j\le n}C_{n+i-j}^{i}\Lambda(j)-C_{n+i}^{i+1}\right)\\ &=\frac{1}{C_{n+i-1}^{i}}\sum_{j\le n}C_{n+i-j}^{i}\Lambda(j)-\frac{n+i}{i+1}\end{aligned} \tag{2.9}$$

Denote $a_{n,j}^{(i)}=C_{n+i-j}^{i}/C_{n+i-1}^{i}$, define $\psi_i(x)=\sum_{j\le x}a_{[x],j}^{(i)}\Lambda(j)$, $\psi_i(x)$ may be viewed as a weighted sum of function $\Lambda(x)$, specially, $\psi_0(x)=\psi(x)$. And (2.9) can be written as

$$\overline{r}^{(i)}(x)=\psi_i(x)-([x]+i)/(i+1). \tag{2.10}$$

Equation (2.10) can be interpreted as

$$\psi_i(x)\sim(x+i)/(i+1), \tag{2.11}$$

and $\overline{r}^{(i)}(x)$ is the error part.

By PC, we have calculated the concrete values of $\overline{r}^{(i)}(n), 1 \le i \le 3$, for $n \le 10^5$, it is somewhat unexpected that $\overline{r}^{(i)}(x)$ are very good at convergence,

| $1 \le \mathrm{n} \le 10^5$ | min | max |
|---|---|---|
| $r(n)$ | -161.501282 | 173.492942 |
| $\overline{r}(n)$ | -5.183956 | 2.717997 |
| $\overline{r}^{(2)}(n)$ | -1.866302 | -0.922313 |
| $\overline{r}^{(3)}(n)$ | -1.428963 | -1.000000 |

This indicates that our estimations (1.4), (1.5) and (2.7) are far from the truth, and likely $\overline{r}^{(i)}(x) \ll x^{\varepsilon}$, and $\lambda_i(x) \ll x^{-1/2+\varepsilon}, 1 \le i \le 3.$

As we seen, the values of $\overline{r}^{(i)}(x)$ are distributed in a very small interval centered at $\overline{c}$, a constant about $-1.2$, so if we retake $x+\overline{c}$ as the dominant of function $\psi(x)$, i.e. $\psi(x) \sim x+\overline{c}$, then new error $r_{\text{new}}(x) = \psi(x) - x - \overline{c} = r(x) - \overline{c}$, so, $\overline{r}_{\text{new}}^{(i)}(x) = \overline{r}^{(i)}(x) - \overline{c}$, and the values of new average errors $\overline{r}_{\text{new}}^{(i)}(x)$ will be re-concentrated at 0-axis.

In view of the observation above that induces us to carry on the difference

$$\overline{r}^{(i)}(n) - \overline{r}^{(i)}(n-1) = \psi_i(n) - \psi_i(n-1) - \frac{1}{i+1} = \sum_{j \le n} (a_{n,j}^{(i)} - a_{n-1,j}^{(i)})\Lambda(j) - \frac{1}{i+1}$$
$$= \frac{1}{i+1}\left( \sum_{j \le n} \frac{C_{n+i-j-1}^{i-1}}{C_{n+i-1}^{i+1}}(j-1)\Lambda(j) - 1 \right)$$

i.e.

$$(i+1)\left(\overline{r}^{(i)}(n) - \overline{r}^{(i)}(n-1)\right) = \sum_{j \le n} \frac{C_{n+i-1-j}^{i-1}}{C_{n+i-1}^{i+1}}(j-1)\Lambda(j) - 1. \qquad (2.12)$$

Denote by $b_{n,j}^{(i)} = (j-1) \cdot C_{n+i-1-j}^{i-1} / C_{n+i-1}^{i+1}$, it is easy to be followed that

$$\sum_{j \le n} b_{n,j}^{(i)} = 1. \qquad (2.13)$$

Define $\hat{\psi}_i(x) = \sum_{j \le x} b_{[x],j}^{(i)} \Lambda(j)$, and write $\hat{r}^{(i)}(x) = (i+1)\left(\overline{r}^{(i)}(x) - \overline{r}^{(i)}(x-1)\right)$, then (2.12) can be written as

$$\hat{r}^{(i)}(x) = \hat{\psi}_i(x) - 1. \qquad (2.14)$$

As $\hat{r}^{(i)}(n)$ is approaching to zero, so (2.14) may be expressed as

$$\hat{\psi}_i(x) \approx 1. \qquad (2.15)$$

Of course, (2.12) can be rewritten as

$$(n-1)\left(\overline{r}^{(i)}(n)-\overline{r}^{(i)}(n-1)\right) = \sum_{j\le n}\frac{C_{n+i-1-j}^{i-1}}{C_{n+i-1}^{i}}(j-1)\Lambda(j)-\frac{n-1}{i+1}, \qquad (2.12')$$

And let $\hat{r}^{(i)}(n)=(n-1)\left(\overline{r}^{(i)}(n)-\overline{r}^{(i)}(n-1)\right)$, $\hat{\psi}_i(n)=\sum_{j\le n}\frac{C_{n+i-1-j}^{i-1}}{C_{n+i-1}^{i}}(j-1)\Lambda(j)$, then (2.12')

becomes

$$\hat{r}^{(i)}(x)=\hat{\psi}_i(x)-\frac{[x]-1}{i+1} \qquad 2.14')$$

And

$$\hat{\psi}_i(x) \sim \frac{x-1}{i+1}. \qquad (2.15')$$

Obviously, $\hat{\psi}_i(n)=\frac{n-1}{i+1}\hat{\psi}_i(n)$, $\hat{r}^{(i)}(n)=\frac{n-1}{i+1}\hat{r}^{(i)}(n)$.

The following are the data by PC for $n \le 10^5$

| $100 \le n \le 10^5$ | min | max |
|---|---|---|
| $\hat{r}^{(1)}(n)$ | -0.089799 | 0.101644 |
| $\hat{r}^{(2)}(n)$ | -0.012375 | 0.007549 |
| $\hat{r}^{(3)}(n)$ | -0.002883 | 0.001493 |
| $\hat{r}^{(4)}(n)$ | -0.001256 | 0.000263 |
| $\hat{r}^{(5)}(n)$ | -0.001183 | 0.000063 |

and

| $1 \le n \le 10^5$ | min | max |
|---|---|---|
| $\hat{r}^{(1)}(n)$ | -159.429591 | 173.815208 |
| $\hat{r}^{(2)}(n)$ | -6.988295 | 8.203225 |
| $\hat{r}^{(3)}(n)$ | -1.520785 | 1.277045 |
| $\hat{r}^{(4)}(n)$ | -0.357921 | 0.278090 |
| $\hat{r}^{(5)}(n)$ | -0.106159 | 0.097080 |

As an interest, in the next, we derive similarly another approximate formula.
At first, we rewrite equation (2.12) as

$$\overline{r}^{(i)}(n)-\overline{r}^{(i)}(n-1)=\frac{1}{n(n-1)}\sum_{j\le n}\frac{P_{n+i-1-j}^{i-1}}{P_{n+i-1}^{i-1}}i(j-1)\Lambda(j)-\frac{1}{(i+1)},$$

i.e.

$$n(n-1)(\overline{r}^{(i)}(n)-\overline{r}^{(i)}(n-1))=\sum_{j\le n}\frac{P_{n+i-1-j}^{i-1}}{P_{n+i-1}^{i-1}}i(j-1)\Lambda(j)-\frac{n(n-1)}{(i+1)}.$$

Let $\mathfrak{r}(n)=n(n-1)(\overline{r}^{(i)}(n)-\overline{r}^{(i)}(n-1))$, then

$$\mathfrak{r}(n)-\mathfrak{r}(n-1)=\sum_{j\le n}\frac{C_{n+i-2-j}^{i-2}}{C_{n+i-1}^{i}}j(j-1)\Lambda(j)-\frac{2(n-1)}{(i+1)}. \tag{2.16}$$

Let $\tilde{r}^{(i)}(n)=\dfrac{\mathfrak{r}(n)-\mathfrak{r}(n-1)}{2}$, denote by $h_{n,j}^{(i)}=\dfrac{C_{n+i-2-j}^{i-2}C_j^2}{C_{n+i-1}^{i}}$, define $\tilde{\psi}_i(x)=\sum_{j\le x}h_{[x],j}^{(i)}\Lambda(j)$, then equation (2.16) will become

$$\tilde{r}^{(i)}(x)=\tilde{\psi}_i(x)-\frac{[x]-1}{(i+1)}. \tag{2.17}$$

In the right-side of (2.16), by $j=n-(n-j)$, separate the sum into two ones, then it is easy to follow that

$$\tilde{r}^{(i)}(n)\ll(\hat{r}^{(i-1)}(n)+\hat{r}^{(i)}(n-1))n.$$

So,

$$\tilde{\psi}_i(x)\sim\frac{x-1}{(i+1)}. \tag{2.18}$$

The following data is by PC for $n\le 10^5$

| $1\le n\le 10^5$ | min | max |
|---|---|---|
| $\tilde{r}^{(2)}(n)$ | -159.856110 | 172.288023 |
| $\tilde{r}^{(3)}(n)$ | -9.331084 | 12.739719 |
| $\tilde{r}^{(4)}(n)$ | -2.853753 | 2.521717 |
| $\tilde{r}^{(5)}(n)$ | -0.856889 | 0.680470 |
| $\tilde{r}^{(6)}(n)$ | -0.299480 | 0.256453 |

## References

1. T.M. Apostol, Introduction to Analytic number theory, Springer-Verlag, 1976
2. L.K. Hua, Introduction to number theory, Sci. Press, 1975, (Chinese)
3. C.T. Pan & C.B. Pan, Introduction to analytic number theory, Sci. Press, 1987, (Chinese)
4. E.C. Titchmarsh, The Theory of the Riemann's zeta-Function, Oxford Clarendon Press, 1986
5. A.Walfisz, Weylsche Exponentialsummen in der Neueren Zahlentheorie, Mathematische Forschungsberichte, XV, V E B Deutscher Verlag der Wissenschaften Berlin, 1963